\documentclass[12pt]{article}
\usepackage[top=2.54cm, bottom=2.54cm, left=2.70cm, right=2.70cm]{geometry}
\usepackage{tikz}
\usepackage[tbtags]{amsmath}
\usepackage{amssymb}
\usepackage{amsthm}
\usepackage{fancyhdr}
\usepackage{latexsym}
\usepackage{mathrsfs}
\usepackage{wasysym}
\usepackage{fancyhdr}
\usepackage{float}
\usepackage{graphicx}
\usepackage[numbers,sort&compress]{natbib}
\usepackage{setspace}
\allowdisplaybreaks[4]

\newtheorem{theorem}{\bf Theorem}[section]
\newtheorem{lemma}[theorem]{Lemma}

\newtheorem{problem}[theorem]{Problem}
\newtheorem{corollary}[theorem]{Corollary}
\newtheorem{conj}[theorem]{Conjecture}

\newtheorem{claim}{\indent Claim}[]
\newtheorem{obser}[theorem]{Observation}

\begin{document}
\begin{spacing}{1.1}
\title{Long antipaths and anticycles in oriented graphs}
\author{Bin Chen$^a$, \, Xinmin Hou$^{a,b,c}$\thanks {Email address: xmhou@ustc.edu.cn(X. Hou)}, \, Xinyu Zhou$^b$\\
\small$^a$Hefei National Laboratory,\\
\small University of Science and Technology of China, Hefei 230088, China\\
\small $^b$School of Mathematical Sciences,\\
\small  University of Science and Technology of China, Hefei 230026, Anhui, China\\
\small$^c$ CAS Key Laboratory of Wu Wen-Tsun Mathematics,\\
\small  University of Science and Technology of China, Hefei 230026, Anhui, China}
\date{}
\maketitle
\noindent{\bf Abstract}\,: Let $\delta^{0}(D)$ be the minimum semi-degree of an oriented graph $D$. Jackson (1981) proved that every oriented graph $D$ with $\delta^{0}(D)\geq k$ contains a directed path of length $2k$ when $|V(D)|>2k+2$,  and a directed Hamilton cycle when $|V(D)|\le 2k+2$.
Stein~(2020) further conjectured that  every oriented graph $D$ with $\delta^{0}(D)>k/2$ contains any  orientated path of length $k$.
Recently, Klimo\u{s}ov\'{a} and Stein (DM, 2023) introduced  the minimum pseudo-semi-degree $\tilde\delta^0(D)$ (a slight weaker than the minimum semi-degree condition as $\tilde\delta^0(D)\ge \delta^0(D))$ and showed that every oriented graph $D$ with  $\tilde\delta^{0}(D)\ge (3k-2)/4$ contains each antipath of length $k$ for $k\geq 3$. 
In this paper, we improve the result of Klimo\u{s}ov\'{a} and Stein by showing that for all $k\geq 2$, every oriented graph with  $\tilde\delta^0(D)\ge(2k+1)/3$ contains either an antipath of length at least $k+1$ or an anticycle of length at least $k+1$. 
Furthermore,  we answer a problem raised by Klimo\u{s}ov\'{a} and Stein in the negative.

{\bf AMS}\,: 05C20; 05C38.


\section{Introduction}

\noindent

All the digraphs considered in this paper are finite and have no loop or parallel arc. 
Let $D=(V, A)$ be a digraph. For an arc $(u, v)\in A(D)$, the vertex $u$ is its \emph{tail} and $v$ is its \emph{head};  we say that $u$ {\em dominates} $v$ or $v$ is {\em dominated} by $u$. 
Write $\delta^{+}(D)$ and $\delta^{-}(D)$ for the \emph{minimum out-degree} and \emph{minimum in-degree} of $D$, respectively. Let $\delta^{0}(D)=\min\{\delta^+(D), \delta^-(D)\}$ be the \emph{minimum semi-degree} of $D$. 

An \emph{oriented graph} is a digraph in which there is at most one arc between any two different vertices. An \emph{orientation} of an undirected graph $G$ is an assignment of exactly one direction to each edge of $G$. A \emph{directed path}\;(resp., \emph{anti-directed path}) is an oriented path such that any two adjacent arcs have the same orientation (resp., opposite orientations). A \emph{directed cycle} (resp., \emph{anti-directed cycle}) is an oriented cycle such that every pair of two adjacent arcs have the same orientation (resp., opposite orientations). An \emph{oriented Hamilton path}\;(resp., \emph{oriented Hamilton cycle}) is an oriented path\;(resp., oriented cycle) passing through all vertices of a given digraph. A digraph $D$ is \emph{strongly connected}\;(\emph{strong} for short) if every two vertices of $D$ are mutually reachable by the directed paths in $D$.

Identifying the sufficient conditions of a digraph for the existence of long oriented paths and oriented cycles is one of the fundamental and well-studied problems in digraph theory. There are many classical and meaningful results on directed Hamilton paths and directed Hamilton cycles. Note that any directed Hamilton cycle contains a directed Hamilton path. In 1960, Ghouila-Houri \cite{G} proved that every strong digraph $D$ on $n$ vertices with $\delta^{+}(D)+\delta^{-}(D)\ge n$ contains a directed Hamilton cycle; in particular,  they showed that each digraph $D$ with $\delta^{0}(D)\geq n/2$ has a directed Hamilton cycle, a directed version of the celebrated Dirac's theorem \cite{D}. In 1972, Woodall \cite{W} verified that for any strong digraph $D$ on $n\geq 2$ vertices, if $d^{+}_{D}(u)+d^{-}(v)\ge n$ for every pair of nonadjacent  vertices $u$ and $v$, then $D$ admitted a directed Hamilton cycle, an analog of the famous Ore's theorem \cite{O} for digraphs.

Thomassen in \cite{T} raised a natural question for oriented graphs: determine the minimum semi-degree that forces a directed Hamilton cycle. H\"{a}ggkvist \cite{H} constructed a family of oriented graphs on $n$ vertices with minimum semi-degree $(3n-5)/8$ containing no directed Hamilton cycle, and conjectured that any oriented graph on $n$ vertices with minimum semi-degree $(3n-4)/8$ has directed Hamilton cycles. In 2009, Keevash, K\"{u}hn and Osthus \cite{KKO09} solved the above conjecture for sufficiently large $n$. In 1993, H\"{a}ggkvist \cite{H} proposed an additional conjecture stating that every oriented graph $D$ on $n$ vertices with the sum of the minimum degree, minimum out-degree and minimum in-degree greater than $(3n-3)/2$ has a Hamilton cycle. In 2008, Kelly, K\"{u}hn and Osthus \cite{KKO08}  approximately verified this conjecture. H\"{a}ggkvist and Thomason in \cite{HT97} conjectured that for any positive $\alpha$, every oriented graph $D$ on $n$ vertices with $\delta^{0}(D)\geq (3/8+\alpha)n$ contains each orientation of a Hamilton cycle, where $n$ is large enough. This was resolved by Kelly \cite{K} in 2011. Particularly, it indicates  for any $k$ such that $3n/4+o(n)\leq k<n$, every oriented graph $D$ on $n$ vertices with $\delta^{0}(D)\geq k/2$ admits any orientation of a path of length $k$.

In 1981, Jackson \cite{J} proved that every oriented graph $D$ with $\delta^{0}(D)\geq k$ contains a directed path of length $2k$ when $|V(D)|>2k+2$,  and a directed Hamilton cycle when $|V(D)|\le 2k+2$.
Stein~\cite{S} further proposed  the following conjecture.
\begin{conj}\textnormal{(\cite{S})} \label{conj}
For each integer $k$, every oriented graph $D$ with $\delta^{0}(D)>k/2$ contains any  oriented path of length $k$.
\end{conj}

In a sense, Conjecture~\ref{conj} is  best possible because of the disjoint union of regular tournaments on $k$ vertices, where $k$ is an odd integer. As shown in \cite{K}, Conjecture~\ref{conj} holds for oriented graphs on $n$ vertices with $k\geq 3n/4+o(n)$.
Jackson \cite{J} showed that Conjecture~\ref{conj} is true for directed paths. Moreover, Fernandes and Lintzmayer \cite{FL} proved that this conjecture is true for all $k\leq 5$. In 2022, Stein and Z\'{a}rate-Guer\'{e}n \cite{SZ} showed that Conjecture~\ref{conj} holds asymptotically for sufficiently large anti-directed paths in sufficiently large oriented graphs.

\begin{theorem}\textnormal{(\cite{SZ})} \label{asympotically}
For all $\alpha\in(0,1)$, there exists $n_{0}$ such that for all $n\geq n_{0}$ and $k\geq \alpha n$, every oriented graph $D$ on $n$ vertices with $\delta^{0}(D)>(1+\alpha)\frac{k}{2}$ contains every anti-directed path of length $k$.
\end{theorem}


In \cite{KS}, Klimo\u{s}ov\'{a} and Stein defined the \emph{minimum pseudo-semi-degree} of a digraph $D$, denoted by $\tilde{\delta}^{0}(D)$, as 
$$\min\{d^+_D(v), d^-_D(v) : v\in V(D)\text{ with }\;d^{+}_{D}(v)>0 \text{ or } d^{-}_{D}(v)>0\}.$$
Clearly, $\delta^0(D)\le\tilde{\delta}^{0}(D)$ and $\tilde{\delta}^{0}(D)=0$ if and only if $A(D)=\emptyset$.  Klimo\u{s}ov\'{a} and Stein~\cite{KS} verified the following result in support of Conjecture~\ref{conj}.

\begin{theorem}\textnormal{(\cite{KS})} \label{KlSt}
Let $k\geq 3$ and $D$ be an oriented graph with $\tilde{\delta}^{0}(D)\geq (3k-2)/4$. Then $D$ contains each anti-directed path of length $k$.
\end{theorem}

An arc-version of Conjecture~\ref{conj} was proposed by Addario-Berry, Havet, Linhares Sales, Reed and Thomass\'{e}~\cite{AHLRT}.

\begin{conj}\textnormal{(\cite{AHLRT})} \label{conj-tree}
For each positive integer $k$, any digraph $D$ with more than $(k-1)|V(D)|$ arcs contains every anti-directed tree on $k+1$ vertices, where an anti-directed tree is an oriented tree in which every vertex has in-degree 0 or out-degree 0.
\end{conj}

As shown in \cite{AHLRT}, Conjecture~\ref{conj-tree} is true for trees with diameter at most 3. In \cite{SZ}, Stein and Z\'{a}rate-Guer\'{e}n approximately settled Conjecture~\ref{conj-tree} for large balanced anti-directed trees in dense oriented graphs. On the other hand, Addario-Berry et al. \cite{AHLRT} showed that every digraph with size greater than $4(k-1)|V(D)|$ contains each anti-directed tree whose largest partition class has at most $k$ vertices. This yields that there must exist each anti-directed path of length $k$ in every digraph $D$ with size larger than $4(\lceil (k+1)/2\rceil-1)|V(D)|$. Klimo\u{s}ov\'{a} and Stein~\cite{KS} improved this result by giving the following theorem.
\begin{theorem}\textnormal{(\cite{KS})} \label{KSsize}
Let $k\geq 1$ and $D$ be an oriented graph with more than $\frac{(3k-4)|V(D)|}2$ arcs. Then $D$ contains each antipath of length $k$.
\end{theorem}

In this article, we further improve Theorems~\ref{KlSt} and~\ref{KSsize} as follows.

\begin{theorem}\label{path-cycle}
Let $D$ be an oriented graph with $\tilde{\delta}^{0}(D)\geq (2k+1)/3$ for any integer $k\geq 2$. Then $D$ contains either an anti-directed path of length at least $k+1$ or an anti-directed cycle of length at least $k+1$.
\end{theorem}

This is an improvement of Theorem~\ref{KlSt}.

 In \cite{SZ}, Klimo\u{s}ov\'{a} and Stein proved the following lemma.
\begin{lemma}\textnormal{(\cite{SZ})}  \label{subdigraph}
	Let $k\geq 1$ and $D$ be a digraph with more than $k|V(D)|$ arcs. Then $D$ contains a subdigraph $D^{*}$ such that $\tilde{\delta}^{0}(D^{*})\geq (k+1)/2$.
\end{lemma}

By Lemma~\ref{subdigraph}, an oriented graph with more than $\frac{(4k-1)|V(D)|}3$ arcs has a subdigraph $D^{*}$ of $D$ such that $\tilde{\delta}^{0}(D^{*})\geq ((4k-1)/3+1)/2=(2k+1)/3$. By Theorem~\ref{path-cycle}, $D^{*}$\;(resp., $D$) contains either an antipath of length at least $k+1$ or an anticycle of length at least $k+1$.  Therefore, we have the following corollary, which is an improvement of Theorem~\ref{KSsize}. 
\begin{corollary}\label{Size}
Let $k\geq 2$ and $D$ be an oriented graph with more than $\frac{(4k-1)|V(D)|}3$ arcs. Then $D$ contains either an anti-directed path of length at least $k+1$ or an anti-directed cycle of length at least $k+1$.
\end{corollary}


Additionally, In \cite{SZ}, Klimo\u{s}ov\'{a} and Stein proposed the following problem.

\begin{problem}\label{counterexample}
Does Theorem~\ref{KlSt} hold for other types of oriented paths except for anti-directed paths?
\end{problem}

We answer this problem in the negative  by the following construction.

\noindent\textbf{Construction of $\mathcal{D}$.} Let $\mathcal{D}$ be a family of oriented bipartite graphs with vertex classes $X$ and $Y$ such that 

(i) $\min\{|X|,|Y|\}\ge (3k-2)/4$ for $k\geq 2$; 

(ii) all arcs are oriented from $X$ to $Y$; 

(iii) $d^{+}_{D}(x)\ge (3k-2)/4$ and $d^{-}_{D}(y)\ge(3k-2)/4$ for any $x\in X$ and $y\in Y$.
\vspace{0.2cm}

Indeed, every oriented graph $D\in \mathcal{D}$ has $\tilde{\delta}^{0}(D)\ge (3k-2)/4$, however,  $D$  contains no directed paths of length two. This implies that $D$ contains no oriented path of length at least two except for the antipaths. Therefore, Construction $\mathcal{D}$ gives a negative answer to Problem~\ref{counterexample}.

The remainder of this paper is organized as follows. In the next section, we introduce additional notation and provide some preliminaries. The proof of Theorem~\ref{path-cycle}  is given in Section 3. We provide some remarks and discussions in the last section.

\section{Preliminaries}






Let $D=(V,A)$ be a digraph. We write $u\rightarrow v$ for $(u,v)\in A(D)$ and $u\nrightarrow v$ for $(u,v)\notin A(D)$.
Let $v\in V(D)$. The \emph{out-neighborhood} of $v$ is the set $N^{+}_{D}(v)=\{u\in V(D):v\rightarrow u\}$, i.e., $d^{+}_{D}(v)=|N^{+}_{D}(v)|$. Similarly, the \emph{in-neighborhood} of $v$ is the set $N^{-}_{D}(v)=\{u\in V(D) : u\rightarrow v\}$, and $d^{-}_{D}(v)=|N^{-}_{D}(v)|$. The vertices in $N^{+}_{D}(v)$ and $N^{-}_{D}(v)$ are called the \emph{out-neighbors} and \emph{in-neighbors} of $v$, respectively. 
For $X, Y\subseteq V(D)$, let $A(X,Y)=\{(x,y)\in A(D) : \, x\in X \text{ and  } y\in Y\}$. As usual, write $D[X]$ for the subdigraph induced by $X$ in $D$.

Recall that an anti-directed path (\emph{antipath} for short) $P$ of length $\ell$ is  a list of $\ell+1$ distinct vertices $v_{0},v_{1},\cdots,v_{\ell}$ such that either $v_{0}\rightarrow v_{1}\leftarrow v_{2}\rightarrow\cdots$ or $v_{0}\leftarrow v_{1}\rightarrow v_{2}\leftarrow\cdots$. It is easily seen that there is only one type antipath of length $\ell$ if $\ell$ is odd, and there are two types antipaths of length $\ell$ if $\ell$ is even.  An anti-directed cycle (\emph{anticycle} for short) of length $\ell$ of $D$ is a list of $\ell$ distinct vertices $u_{1},u_{2},\cdots,u_{\ell}$ such that either $u_{i}\rightarrow \{u_{i-1},u_{i+1}\}$ for odd $i$ or $u_{i}\rightarrow \{u_{i-1},u_{i+1}\}$ for even $i$, where $i\in\{1,2,\cdots,\ell\}$. Note that the length of every anticycle is even, i.e., there is only one type of anticycle.


\noindent


The following observation is very useful.
\begin{obser}\label{longest antipaths}
 Let $D$ be an oriented graph and $P=x_{1}x_{2}\cdots x_{k}$ be a longest antipath in $D$ for $k\geq 1$. Then we have

 $(i)$. If $x_{1}\rightarrow x_{2}$ in $P$, then $N^{+}_{D}(x_{1})\subseteq V(P)$;

 $(ii)$. If $x_{2}\rightarrow x_{1}$ in $P$, then $N^{-}_{D}(x_{1})\subseteq V(P)$.
\end{obser}

The next lemma is weaker than Theorem~\ref{path-cycle}; however, the technique used here reveals the main idea of the proof of Theorem~\ref{path-cycle} and the result will be used  at the beginning of that proof.

\begin{lemma}\label{k+1}
	Let $D$ be an oriented graph with $\tilde{\delta}^{0}(D)\geq k$ for $k\geq 2$. Then $D$ contains either an antipath of length at least $k+1$ or an anticycle of length at least $k+1$.
\end{lemma}
\noindent\textbf{Proof.}
For the sake of contradiction, suppose that there is neither an antipath of length at least $k+1$ nor an anticycle of length at least $k+1$. Note that $A(D)\neq \emptyset$ since $\tilde{\delta}^{0}(D)\geq k\geq 2$. Let $P=x_{0}x_{1}\cdots x_{\ell}$ be a longest antipath in $D$. Then $\ell\leq k$. Assume, without loss of generality, that $x_{0}\rightarrow x_{1}$. By Observation~\ref{longest antipaths}, we have $N^{+}_{D}(x_{0})\subseteq V(P)$. As $x_{0}\rightarrow x_{1}$, then $d^{+}_{D}(x_{0})\geq \tilde{\delta}^{0}(D)\geq k$, implying that $\ell \geq k$. Together with $\ell\leq k$, we obtain $\ell=k$. 
This forces that $x_0\rightarrow x_k$.

If $k$ is odd, then $x_{k-1}\rightarrow x_{k}$.
Thus $C=x_{0}x_{1}\cdots x_{k}x_{0}$ is an anticycle of length $k+1$, a contradiction. 

If $k$ is even, then $x_{k}\rightarrow x_{k-1}$. By Observation~\ref{longest antipaths}, we also have $N^{+}_{D}(x_{k})\subseteq V(P)$. By $d^{+}_{D}(x_{k})\geq \tilde{\delta}^{0}(D)\geq k$, we have $x_{k}\rightarrow x_{0}$. Thus $x_0\nrightarrow x_k$, which is also a contradiction.  
\hfill $\blacksquare$

\section{The proof of Theorem~\ref{path-cycle}}




\noindent



One can check directly from Lemma~\ref{k+1} that the result holds for $k=2$ or $3$ because $\tilde{\delta}^{0}(D)$ is an integer. Now assume $k\ge 4$. Suppose to the contrary that there is an oriented graph $D$ with $\tilde{\delta}^0(D)\geq (2k+1)/3$, however, $D$ does not contain any antipath  or anticycle of length at least $k+1$.
Let $$\mathscr{P}=\{P:P\;\text{is}\;\text{a}\;\text{longest}\;\text{antipath}\;\text{in}\;D\}$$ and $$\mathscr{C}=\{C:C\;\text{is}\;\text{a}\;\text{longest}\;\text{anticycle}\;\text{in}\;D\}.$$
Since $k\geq 4$ and $\tilde{\delta}^{0}(D)\geq (2k+1)/3$, we have $A(D)\neq \emptyset$. Therefore, $\mathscr{P}\neq \emptyset$. Let $t$ be the length of a longest antipath in $D$. Then $t\leq k$ according to the assumption.

We first show that $\mathscr{C}\neq \emptyset$.

\begin{claim}\label{anticycle}
$\mathscr{C}\neq \emptyset$.
\end{claim}
\begin{proof}

It is sufficient to verify that $D$ contains anticycles. Choose a longest antipath $P=x_{0}x_{1}\cdots x_{t}$ in $D$. By symmetry, we may assume that $x_{0}\rightarrow x_{1}$.

If $t$ is odd, then $x_{i}\rightarrow \{x_{i-1},x_{i+1}\}$ for all $i\in \{2,4,\cdots, t-1\}$. Since $x_{0}\rightarrow x_{1}$, by Observation~\ref{longest antipaths}, we have $N^{+}_{D}(x_{0})\subseteq V(P)$.  If there is a $x_{i}\in \{x_{3},x_{5},\cdots,x_{t}\}$ such that $x_{0}\rightarrow x_{i}$, then  $x_{0}x_{1}\cdots x_{i}x_{0}$ is an anticycle. Therefore, $N^{+}_{D}(x_{0})\subseteq \{x_{1},x_{2},x_{4},\cdots, x_{t-1}\}$. Hence $d^{+}_{D}(x_{0})\leq 1+(t-1)/2\leq (k+1)/2$. As $\tilde{\delta}^{0}(D)\geq (2k+1)/3$ and $d^{+}_{D}(x_{0})\geq 1$, we have $(2k+1)/3\le d^{+}_{D}(x_{0})\le (k+1)/2$, this is impossible since $k\ge 4$.

Now assume  $t$ is even.
Then, apparently, $\{x_{i-1},x_{i+1}\}\rightarrow x_{i}$ for all $i\in \{1,3,\cdots,t-1\}$. Again according to Observation~\ref{longest antipaths}, we have $N^{+}_{D}(x_{0})\subseteq V(P)$ and $N^{+}_{D}(x_{t})\subseteq V(P)$. If $N^{+}_{D}(x_{0})\cap \{x_{3},x_{5},\cdots,x_{t-1}\}\not=\emptyset$ or $N^{+}_{D}(x_{t})\cap \{x_{1},x_{3},\cdots,x_{t-3}\}\not=\emptyset$, then we have an anticycle in $D$ as in the above case. Note that there exists at most one arc between $x_{0}$ and $x_{t}$. Without loss of generality, we may assume that $x_{0}\nrightarrow x_{t}$. Hence $N^{+}_{D}(x_{0})\subseteq \{x_{1},x_{2},x_{4},\cdots,x_{t-2}\}$, which yields that $d^{+}_{D}(x_{0})\leq 1+(t-2)/2\leq k/2$. Together with $d^{+}_{D}(x_{0})\geq (2k+1)/3$, we have a contradiction again.
\end{proof}

Let $s$ be the length of a longest anticycle in $D$. Then, $s$ must be even, and $s\leq k$.

\begin{claim}\label{length}
$t\geq s$.
\end{claim}
\begin{proof}

It suffices to find an antipath of length at least $s$. Suppose that there is no such antipath. Choose a longest anticycle $C=y_{1}y_{2}\cdots y_{s}y_{1}$ in $D$ such that $y_{i}\rightarrow \{y_{i-1}, y_{i+1}\}$ for all $i\in \{1,3,\cdots,s-1\}$, where the sum in the index is under module $s$. If there exists $y_{i}$ with $i\in \{1,3,\cdots,s-1\}$ such that $y_{i}\rightarrow v$ for some $v\in V(D)\backslash V(C)$, then $y_{i+1}y_{i+2}\cdots y_{s}y_{1}\cdots y_{i}v$ is an antipath of length $s$, a contradiction.
Hence $N^{+}_{D}(y_{i})\subseteq V(C)$ for all $i\in \{1,3,\cdots,s-1\}$. Analogously, we have $N^{-}_{D}(y_{j})\subseteq V(C)$ for all $j\in \{2,4,\cdots,s\}$. Denote by $U=\{y_{i}\neq y_{s}:y_{i+2}\in N^{+}_{D}(y_{1})\}$. Since $d^{+}_{D}(y_{1})\geq (2k+1)/3$ and $y_{1}\rightarrow y_{s}$, we have $|U|\geq (2k-2)/3$.
We claim that $N^{-}_{D}(y_{s})\cap U\neq \emptyset$. Suppose that $N^{-}_{D}(y_{s})\cap U=\emptyset$.
Due to $N^{-}_{D}(y_{s})\subseteq V(C)$ and $|V(C)|=s\leq k$, we have $d^{-}_{D}(y_{s})\leq s-1-|U|\leq (k+2)/3$.
However, $d^{-}_{D}(y_{s})\ge \tilde{\delta}^{0}(D)\geq (2k+2)/3$ because $y_{1}\rightarrow y_{s}$, a contradiction.
Therefore, $N^{-}_{D}(y_{s})\cap U\neq \emptyset$. Namely, there exist $y_{i},y_{i-2}$ such that $y_{1}\rightarrow y_{i}$ and $y_{i-2}\rightarrow y_{s}$, where $i\neq 2$.

If $i$ is odd, then  $y_{i}\rightarrow \{y_{i-1},y_{i+1}\}$ on $C$. Hence $d^{+}_{D}(y_{i})\ge \tilde{\delta}^{0}(D)\geq (2k+1)/3$. Since $y_1\rightarrow y_i$, we have $d^{-}_{D}(y_{i})\ge \tilde{\delta}^{0}(D)\geq (2k+1)/3$. As $N^{+}_{D}(y_{i})\subseteq V(C)$, we obtain $|N^{-}_{D}(y_{i})\cap V(C)|\leq s-1-(2k+1)/3\leq (k-4)/3$.
Thus we have $|N^{-}_{D}(y_{i})\cap (V(D)\backslash V(C))|\geq (2k+1)/3-(k-4)/3=(k+5)/3>0$. Choose $v\in N^{-}_{D}(y_{i})\cap (V(D)\backslash V(C))$. Then $d^{+}_{D}(v)\geq (2k+1)/3$ as $v\rightarrow y_{i}$. Moreover, $v\nrightarrow y_{j}$ for any $j\in\{2,4,\cdots,s\}$ since, otherwise, we will obtain an antipath $vy_jy_{j-1}\dots y_1y_sy_{s-1}\dots y_{j+1}$ of length at least $s$. It follows that $|N^{+}_{D}(v)\cap V(C)|\subseteq \{y_{1},y_{3},\cdots,y_{s-1}\}$. This indicates that $|N^{+}_{D}(v)\backslash V(C)|\geq (2k+1)/3-s/2\geq (k+2)/6>0$.
Pick $u\in N^{+}_{D}(v)\backslash V(C)$. Then $uvy_{i}y_{1}\cdots y_{i-2}y_{s}y_{s-1}\cdots y_{i+1}$ is an antipath of length $s$, a contradiction.

If $i$ is even, then $y_{i-2}\rightarrow y_{s}$ and $\{y_{i-3},y_{i-1}\}\rightarrow y_{i-2}$. We also have $d^{+}_{D}(y_{i-2})\geq (2k+1)/3$ and $d^{-}_{D}(y_{i-2})\geq (2k+1)/3$. In addition, $|N^{-}_{D}(y_{i-2})\cap V(C)|\geq (2k+1)/3$ implies that $|N^{+}_{D}(y_{i-2})\cap V(C)|\leq s-1-(2k+1)/3\leq (k-4)/3$. Therefore, $|N^{+}_{D}(y_{i-2})\cap (V(D)\backslash V(C))|\geq (k+5)/3> 0$. Let $v\in V(D)\backslash V(C)$ such that $y_{i-2}\rightarrow v$. Then  $d^{-}_{D}(v)\geq (2k+1)/3$ and $y_{j}\nrightarrow v$ for all $j\in \{1,3,\cdots,s-1\}$. It follows that $|N^{-}_{D}(v)\cap (V(D)\backslash V(C))|\geq (2k+1)/3-s/2\geq (k+2)/6>0$. As a consequence, we can find a vertex $u\in V(D)\backslash V(C)$ such that $u\neq v$ and $u\rightarrow v$. One can see that $uvy_{i-2}y_{s}y_{s-1}\cdots y_{i}y_{1}\cdots y_{i-3}$ is an antipath of length $s$, a contradiction too.
\end{proof}

\begin{claim}\label{odd length}
$t$ is odd.
\end{claim}
\begin{proof}
Suppose to the contrary that $t$ is even. 
Choose a longest antipath $Q=z_{0}z_{1}\cdots z_{t}$. By symmetry, we may assume $z_{0}\rightarrow z_{1}$. Note that $s\leq t\leq k$ and $\{z_{i-1},z_{i+1}\}\rightarrow z_{i}$ for all $i\in\{1,3,\cdots,t-1\}$. As $Q$ is a longest antipath in $D$ and $t$ is even, according to Observation~\ref{longest antipaths}, we deduce that $N^{+}_{D}(z_{0})\subseteq V(Q)$ and $N^{+}_{D}(z_{t})\subseteq V(Q)$. We assert that $z_{0}\nrightarrow z_{t}$. If not, we have $z_{0}\rightarrow z_{t}$. It is apparent that $Q^{*}=z_{t}z_{0}z_{1}\cdots z_{t-1}$ is an antipath of length $t$, that is, $Q^{*}$ is also a longest antipath in $D$. As $z_{0}\rightarrow z_{t}$, again according to Observation~\ref{longest antipaths}, one can see that $N^{-}_{D}(z_{t})\subseteq V(Q^{*})=V(Q)$. This implies that $N^{+}_{D}(z_{t})\cup N^{-}_{D}(z_{t})\subseteq V(Q)$. Due to $d^{+}_{D}(z_{t})\geq (2k+1)/3$ and $d^{-}_{D}(z_{t})\geq (2k+1)/3$, we obtain $(4k+2)/3\leq t\leq k$. This contradicts to $k\geq 4$. We thus deduce that $z_{0}\nrightarrow z_{t}$. Similarly, one can also show that $z_{t}\nrightarrow z_{0}$.

Let $W=\{z_{i-1}\neq z_{0}: z_{i}\in N^{+}_{D}(z_{0})\}$. We claim that $N^{+}_{D}(z_{t})\cap W\neq \emptyset$. Conversely, suppose $N^{+}_{D}(z_{t})\cap W= \emptyset$. Then $d^{+}_{D}(z_{t})\leq t-(2k+1)/3\leq (k-1)/3<(2k+1)/3$ as $k\geq 4$, a contradiction. This indicates that for some $i\in \{2,3,\cdots,t-1\}$, $z_{0} \rightarrow z_{i}$ and $z_{t}\rightarrow z_{i-1}$. If $i$ is odd, then both $z_{i-1}z_{i-2}\cdots z_{0}z_{i}z_{i+1}\cdots z_{t}$ and $z_{i-2}z_{i-3}\cdots z_{0}z_{i}z_{i+1}\cdots z_{t}z_{i-1}$ are longest antipaths in $D$. Therefore, we have $N^{+}_{D}(z_{i-1})\cup N^{-}_{D}(z_{i-1})\subseteq V(Q)$ by Observation~\ref{longest antipaths}. It follows that $(4k+2)/3\leq t\leq k$, which is impossible as $k\geq 4$. If $i$ is even, then both $z_{i}z_{i+1}\cdots z_{t}z_{i-1}z_{i-2}\cdots z_{0}$ and $z_{i+1}z_{i+2}\cdots z_{t}z_{i-1}z_{i-2}\cdots z_{0}z_{i}$ are longest antipaths in $D$. We obtain $N^{+}_{D}(z_{i})\cup N^{-}_{D}(z_{i})\subseteq V(Q)$ via Observation~\ref{longest antipaths}. Hence, $(4k+2)/3\leq t\leq k$, a contradiction too.
\end{proof}

In the remainder of the proof, we let $f(k)=(2k+1)/3$.

Let $R=r_{0}r_{1}\cdots r_{t}$ be a longest antipath in $D$. Note first that $t$ is odd according to Claim~\ref{odd length} and $t\leq k$.
By symmetry, we may assume that $r_{0}\rightarrow r_{1}$. Thus we have $r_{i}\rightarrow \{r_{i-1},r_{i+1}\}$ for each $i\in \{2,4,\cdots,t-1\}$. As $R$ is a longest antipath with $r_{0}\rightarrow r_{1}$ and $r_{t-1}\rightarrow r_{t}$, we deduce that $N^{+}_{D}(r_{0})\subseteq V(R)$ and $N^{-}_{D}(r_{t})\subseteq V(R)$ by Observation~\ref{longest antipaths}. In addition, $d^{+}_{D}(r_{0})\geq \tilde{\delta}^{0}(D)\geq f(k)$ and $d^{-}_{D}(r_{t})\geq \tilde{\delta}^{0}(D)\geq f(k)$.

If $r_{0}\rightarrow r_{t}$, then $r_{0}r_{1}r_{2}\cdots r_{t}r_{0}$ is an anticycle of length $t+1$, a contradiction to our assuming  that any longest anticycle in $D$ has length $s\leq t$. Therefore, $r_{0}\nrightarrow r_{t}$.


Denote by $X_{1}=\{r_{1},r_{3},\cdots,r_{t}\}$ and $X_{2}=\{r_{0},r_{2},\cdots,r_{t-1}\}$. Let $Y_{1}=N^{+}_{D}(r_{0})\cap X_{1}$ and $Y_{2}=\{r_{i}:r_{i+1}\in Y_{1}\}$. Then $|Y_{1}|=|Y_{2}|$ and $Y_{2}\subseteq X_2$. For every $r_i\in Y_2$, it is easy to check that $r_{i}r_{i-1}\cdots r_{0}r_{i+1}r_{i+2}\cdots r_{t}$ is a longest antipath in $D$ with  $r_{i}\rightarrow r_{i-1}$. Thus every $r_{i}\in Y_{2}$ satisfies  $N^{+}_{D}(r_{i})\subseteq V(R)$. Recall that $d^{+}_{D}(r_{0})\geq f(k)$ and $N^{+}_{D}(r_{0})\subseteq V(R)$. Thus
$$|Y_{1}|\geq d^{+}_{D}(r_{0})-(|X_{2}|-1)\geq f(k)-(t-1)/2\geq f(k)-(k-1)/2.$$
This  yields $|Y_{2}|\geq f(k)-(k-1)/2$. Note that $D[Y_{2}]$ is an oriented subgraph of $D$. Thus $\delta^{+}(D[Y_{2}])\leq (|Y_{2}|-1)/2$. Hence we can choose  some vertex $r_{j}\in Y_{2}$ such that $|N^{+}_{D}(r_{j})\cap Y_{2}|\leq (|Y_{2}|-1)/2$.
Therefore, we have
\begin{align}\label{EQ:1}
|N^{+}_{D}(r_{j})\cap X_{2}|
&\leq|N^{+}_{D}(r_{j})\cap Y_{2}|+|X_{2}\backslash Y_{2}|\nonumber\\
&\leq(|Y_{2}|-1)/2+|X_{2}|-|Y_{2}|\nonumber\\
&=|X_{2}|-(|Y_{2}|+1)/2\nonumber\\
&=(t-|Y_{2}|)/2\\
&\leq \frac12\left(\frac{3k-1}{2}-f(k)\right).\nonumber
\end{align}
Note again that $N^{+}_{D}(r_{j})\subseteq V(R)$ and $d^{+}_{D}(r_{j})\geq \tilde{\delta}^{0}(D)\geq f(k)$. As $V(R)=X_{1}\cup X_{2}$, we obtain that
\begin{equation}\label{EQ: e1}
|N^{+}_{D}(r_{j})\cap X_{1}|= d^{+}_{D}(r_{j})-|N^{+}_{D}(r_{j})\cap X_{2}|\geq \frac 12\left(3f(k)-\frac{3k-1}{2}\right).
\end{equation}
For convenience, write $R^{1}=r_{0}^{1}r_{1}^{1}r_{2}^{1}\cdots r_{t}^{1}$ for $r_{j}r_{j-1}\cdots r_{0}r_{j+1}r_{j+2}\cdots r_{t}$. Denote by $X_{1}^{1}=\{r_{1}^{1},r_{3}^{1},\cdots,r_{t}^{1}\}$ and $X_{2}^{1}=\{r_{0}^{1},r_{2}^{1},\cdots,r_{t-1}^{1}\}$. It is clear that $X_{i}^{1}=X_{i}$ for $i=1,2$, and $r_{i}^{1}\rightarrow \{r_{i-1}^{1},r_{i+1}^{1}\}$ for every $i\in \{2,4,\cdots,t-1\}$. By (\ref{EQ: e1}), we also have $|N^{+}_{D}(r_{0}^{1})\cap X_{1}^{1}|\geq (3f(k)-\frac{3k-1}{2})/2$.
As $r_{0}^{1}r_{1}^{1}r_{2}^{1}\cdots r_{t}^{1}$ is a longest antipath and $t$ is odd, we have $r^{1}_{0}\nrightarrow r^{1}_{t}$ too. Let $Y_{1}^{1}=N^{+}_{D}(r_{0}^{1})\cap X_{1}^{1}$ and $Y_{2}^{1}=\{r_{i}^{1}:r_{i+1}^{1}\in Y_{1}^{1}\}$. Then $|Y_{1}^{1}|=|Y_{2}^{1}|\geq (3f(k)-\frac{3k-1}{2})/2$. For the same reason as above, there exists a vertex $r_{j}^{1}\in Y_{2}^{1}$ such that $|N^{+}_{D}(r_{j}^{1})\cap Y_{2}^{1}|\leq (|Y_{2}^{1}|-1)/2$. 
As in (\ref{EQ:1}), we have
$$|N^{+}_{D}(r_{j}^1)\cap X_{2}^{1}|\leq (t-|Y_2^1|)/2\le \frac 14\left(\frac{7k-1}{2}-3f(k)\right),$$
and $$|N^{+}_{D}(r_{j}^{1})\cap X_{1}^{1}|\geq \frac 14\left(7f(k)-\frac{7k-1}{2}\right).$$
Set $r_{0}^{2}=r_{j}^{1}$.
Repeating the previous process $\alpha$ times, where $\alpha=\lceil \log_{2}k\rceil$. We obtain a longest antipath $R^{\alpha}=r_{0}^{\alpha}r_{1}^{\alpha}r_{2}^{\alpha}\cdots r_{t}^{\alpha}$. 
Denote by $X_{1}^{\alpha}=\{r_{1}^{\alpha},r_{3}^{\alpha},\cdots,r_{t}^{\alpha}\}$ and $X_{2}^{\alpha}=\{r_{0}^{\alpha},r_{2}^{\alpha},\cdots,r_{t-1}^{\alpha}\}$. Obviously, $X_{i}^{\alpha}=X_{i}$ for $i=1,2$, and $r_{i}^{\alpha}\rightarrow \{r_{i-1}^{\alpha},r_{i+1}^{\alpha}\}$ for every $i\in \{2,4,\cdots,t-1\}$. Moreover, we have$$|N^{+}_{D}(r_{0}^{\alpha})\cap X_{1}^{\alpha}|\geq\frac 1{2^{\alpha}}\left((2^{\alpha+1}-1)f(k)-\frac{(2^{\alpha+1}-1)k-1}{2}\right).$$
Let $$g(k)=f(k)+\frac1{2^{\alpha}}\left((2^{\alpha+1}-1)f(k)-\frac{(2^{\alpha+1}-1)k-1}{2}\right).$$ 
By combining $f(k)=(2k+1)/3$ and $\alpha=\lceil \log_{2}k\rceil$, we have 
\begin{align*}
g(k)
&=\frac1{2^{\alpha}}\left((3*2^{\alpha}-1)f(k)-\frac{(2^{\alpha+1}-1)k-1}{2}\right)\\
&=\frac1{2^{\alpha}}\left(2^{\alpha}(2k+1)-\frac{(2k+1)}3-2^{\alpha}k+\frac {(k+1)}2\right)\\
&=k+1-\frac{(k-1)}{6\cdot 2^{\alpha}}\\
&>k.
\end{align*}
Now let $\widetilde{X}=\{r_{i}^{\alpha}\neq r_{0}^{\alpha}:r_{i+1}^{\alpha}\in N^{+}_{D}(r_{0}^{\alpha})\cap X_{1}^{\alpha}\}$. Then $|\widetilde{X}|=|N^{+}_{D}(r_{0}^{\alpha})\cap X_{1}^{\alpha}|-1$. We claim that $N^{-}_{D}(r_{t}^{\alpha})\cap  \widetilde{X}\neq\emptyset$.  
Otherwise, we have $$d_D^-(r_t^\alpha)\le t-1-|\widetilde{X}|\le k-1-(|N^{+}_{D}(r_{0}^{\alpha})\cap X_{1}^{\alpha}|-1)=k-|N^{+}_{D}(r_{0}^{\alpha})\cap X_{1}^{\alpha}|.$$ 
This implies that $k\ge d_D^-(r_t^\alpha)+|N^{+}_{D}(r_{0}^{\alpha})\cap X_{1}^{\alpha}|\ge g(k)>k$, a contradiction.
Choose $r_{i-1}^{\alpha}\in N^{-}_{D}(r_{t}^{\alpha})\cap  \widetilde{X}$ for some $i\in \{3,5,\cdots,t-2\}$. Then we have $r_{0}^{\alpha}\rightarrow r_{i}^{\alpha}$ and $r_{i-1}^{\alpha}\rightarrow r_{t}^{\alpha}$. Therefore, we have an anticycle $r_{0}^{\alpha}r_{1}^{\alpha}\cdots r_{i-1}^{\alpha}r_{t}^{\alpha}r_{t-1}^{\alpha}\cdots r_{i}^{\alpha}r_{0}^{\alpha}$ of length $t+1$,  a contradiction to Claim~\ref{length}. 

This completes our proof of the theorem.       \hfill $\blacksquare$
\vspace{0.2cm}

\section{Remarks and Discussions}
By constructing family $\mathcal{D}$, we know that a high minimum pseudo-semi-degree does not guarantee long  paths of any orientation. However, it is  still interesting to address the following problem. 
\begin{problem}\label{p1}
	For each integer $k\ge 2$, every oriented graph $D$ with $\tilde\delta^{0}(D)>k/2$ contains an anti-directed path of length $k$.
\end{problem}

\subsection*{Acknowledgements}

\noindent
This work was supported by the National Key Research and Development Program of China (2023YFA1010200), the National Natural Science Foundation of China (No. 12071453),  and the Innovation Program for Quantum Science and Technology (2021ZD0302902).

\vskip 3mm
\end{spacing}

\vspace{0.3cm}









\begin{thebibliography}{99}

\bibitem{AHLRT} L. Addario-Berry, F. Havet, C. Linhares Sales, B. Reed and S. Thomass\'{e}, Oriented trees in digraphs, Discrete Math. 313 (2013) 967-974


\bibitem{BJG} J. Bang-Jensen, G. Gutin, Digraphs: Theory, Algorithms and Applications, second ed., Springer, London, 2009.

\bibitem {D} G.A. Dirac, Some theorems on abstract graphs, Proc. Lond. Math. Soc. 2(1) (1952) 69-81

\bibitem {FL} C.G. Fernandes and C.N. Lintzmayer, Personal communication

\bibitem {G} A. Ghouila-Houri, Une condition suffisante d'existence d'un circuit hamiltonien, C. R. Math. Acad. Sci. Paris. 25 (1960) 495-497

\bibitem {H} R. H\"{a}ggkvist, Hamilton cycles in oriented graphs, Comb. Probab. Comput. 2 (1993) 25-32




\bibitem {HT97} R. H\"{a}ggkvist and A. Thomason, Oriented Hamilton cycles in oriented graphs, in: Combinatorics, Geometry and Probability, Cambridge University Press. (1997) 339-353


\bibitem {J} B. Jackson, Long paths and cycles in oriented graphs, J. Graph Theory. 5(2) (1981) 145-157


\bibitem {KKO09} P. Keevash, D. K\"{u}hn and D. Osthus, An exact minimum degree condition for Hamilton cycles in oriented graphs, J. Lond. Math. Soc. 79 (2009) 144-166



\bibitem {K} L. Kelly, Arbitrary orientations of Hamilton cycles in oriented graphs, Electron. J. Comb. 18(1) (2011) P186

\bibitem {KKO08} L. Kelly, D. K\"{u}hn and D. Osthus, A Dirac type result on Hamilton cycles in oriented graphs, Comb. Probab. Comput. 17 (2008) 689-709


\bibitem {KS} T. Klimo\u{s}ov\'{a} and M. Stein, Antipaths in oriented graphs, Discrete Math. 346 (2023) 113515

\bibitem {O} O. Ore, Note on Hamilton circuits, Amer. Math. Monthly. 67 (1960) 55

\bibitem {S} M. Stein, Tree containment and degree conditions, in: A.M. Raigorodskii, M.T. Rassias (Eds.), Discrete Mathematics and Applications, in: Springer Optimization and Its Applications, Springer. 165 (2020) 459-486




\bibitem {SZ} M. Stein and C. Z\'{a}rate-Guer\'{e}n, Antidirected subgraphs of oriented graphs, arXiv:2212.00769v1 (2022)

\bibitem {T}  C. Thomassen, Long cycles in digraphs with constraints on the degrees, in: Surveys in Combinatorics, in: London Math. Soc. Lecture Notes, Cambridge University Press 38 (1979) 211-228


\bibitem {W} D. Woodall, Sufficient conditions for cycles in digraphs, Proc. Lond. Math. Soc. 24 (1972) 739-755

\end{thebibliography}
\end{document}